\documentclass{amsart}

\def \Re{\mbox{${\mathcal Re}$}}
\def \r{\mbox{${\mathbb R}$}}
\def \C{\mbox{${\mathbb C}$}}
\def \h{\mbox{${\mathbb H}$}}

\newtheorem{thm}{Theorem}[section]

\newtheorem{lem}[thm]{Lemma}

\theoremstyle{remark}
\newtheorem{rem}[thm]{Remark}

\theoremstyle{definition}

\newtheorem{ex}[thm]{Example}

\numberwithin{equation}{section} \numberwithin{thm}{section}

\begin{document}


\title{On the Bj\"{o}rling problem in a three-dimensional Lie group}

\author{Francesco Mercuri}
\address{Departamento de Matem\'{a}tica, C.P. 6065\\
IMECC, UNICAMP, 13081-970, Campinas, SP\\ Brasil}

\email{mercuri@ime.unicamp.br}

\author{Irene I. Onnis}

\address{Departamento de Matem\'{a}tica, C.P. 668\\ ICMC,
USP, 13560-970, S\~{a}o Carlos, SP\\ Brasil}

\email{onnis@icmc.usp.br}

\thanks{Work partially supported by RAS, INdAM, FAPESP and CNPq}

\begin{abstract} We prove existence and uniqueness of the solution of the Bj\"{o}rling problem
for minimal surfaces in a three-dimensional Lie group.
\end{abstract}

\maketitle

\section{Introduction}

The  Weierstrass representation formula for minimal surfaces in
$\r^3$ has been a fundamental tool for producing examples and
proving general properties of such surfaces, since the surfaces can be
parametrized by holomorphic data. In \cite{MMP} the authors describe
a general Weierstrass representation formula for minimal surfaces in
an arbitrary Riemannian manifold. The partial differential equations
involved are, in general, too complicated to be solved explicitly.
However, for particular ambient manifolds, such as the Heisenberg
group, the hyperbolic space and the product of the hyperbolic plane
with $\r$, the equations are more workable and the formula can be
used to produce examples (see \cite{K}, \cite{MMP}).

In this note we will show how this formula can be used, at least if
the ambient manifold is a $3$-dimensional Lie group, in order to
prove existence and uniqueness of the solution of the Bj\"{o}rling
problem. We also give some examples for the case in which the
ambient manifold is the Heisenberg group $\h_3$\ or $\h^2 \times
\r$, the product of the hyperbolic plane and the real line.

We thank the referee for his appropriate and useful observations. In
particular the last section is inspired by his comments.

\section{The Weierstrass representation formula}

The arguments will be essentially local so we will consider, as
ambient manifold $M$, the space $\r^3$ with a Riemannian metric $g =
(g_{ij})$. We will denote by $\Omega \subseteq \C \cong \r^2$ a {\em
simply connected} domain with a complex coordinate $z = u+iv$, $u,v
\in \r$, and by:
$$\frac{\partial}{\partial z} := \frac{1}{2}\ \Big(\frac{\partial}{\partial u} -
i\frac{\partial}{\partial v}\Big);\;\qquad
\frac{\partial}{\partial \overline{z}} :=
\frac{1}{2}\Big(\frac{\partial}{\partial u} +
i\frac{\partial}{\partial v}\Big)
$$
the complex derivatives.

In this situation, the general Weierstrass representation formula
can be stated as follows:
\begin{thm} {\em (see \cite{MMP} for a proof)}\label{wg}
Let $f: \Omega \rightarrow M$ be a conformal minimal immersion and
$g=(g_{ij})$ be the induced metric. The complex tangent vector:
$$\frac{\partial f}{\partial z} := \phi := \sum_i \phi_i\frac{
\partial}{\partial x_i},\qquad \phi_i:\Omega\rightarrow \C,$$ has the
following properties:
\begin{enumerate}
\item $\sum_{i,j} g_{ij}\phi_i\overline{\phi_j} \neq 0,
$
\item $\sum_{i,j} g_{ij}\phi_i{\phi_j} = 0,$
\item $\displaystyle{\frac{\partial \phi_i}{\partial \overline{z}} + \sum_{j,k}
\Gamma_{jk}^i\overline{\phi_j}\phi_k = 0}$,
\end{enumerate}
where $\{\Gamma_{jk}^i\}$ are the Christoffel symbols of the
Riemannian connection.

Conversely, given functions $\phi_i: \Omega \rightarrow\C$ that
verify the above conditions, then the map:
$$f: \Omega \rightarrow M,\qquad f_i(z) = 2\,\Re\int_{z_0}^z \phi_i\, dz,$$
is a well defined conformal minimal immersion of $\Omega$ into $M$
(here $z_0$ is an arbitrary fixed point of $\Omega$ and the integral
is along any curve joining $z_0$\ to $z$).
\end{thm}

\begin{rem} The first condition of Theorem~\ref{wg} tells us that $f$ is an immersion,
the second that $f$ is conformal and the last one that $f$ is
minimal. The last condition is called the {\em holomorphicity
condition} since it is the local coordinates version of the
condition: $\tilde{\nabla}_{\frac{\partial}{\partial \bar{z}}}\phi
=0$, where $\tilde{\nabla}$ is the induced connection on the
pull-back bundle $f^{*}(TM\otimes\C)$. In fact, we have that the
section $\phi$ is holomorphic if and only if
\begin{equation}\begin{aligned}
\tilde{\nabla}_{\frac{\partial}{\partial \bar{z}}}\Big(\sum_i \phi_i
\frac{\partial}{\partial x_i}\Big)&=\sum_i\Big\{\frac{\partial
\phi_i}{\partial \bar{z}}\frac{\partial}{\partial
x_i}+\phi_i\nabla_{\frac{\partial f}{\partial
\bar{z}}}\frac{\partial}{\partial x_i}\Big\}\\&=
\sum_i\Big\{\frac{\partial \phi_i}{\partial
\bar{z}}\frac{\partial}{\partial x_i}+ \phi_i\nabla_{\sum_j
\overline{\phi_j}\frac{\partial }{\partial
x_j}}\frac{\partial}{\partial x_i}\Big\}\\&=
\sum_i\Big\{\frac{\partial \phi_i}{\partial
\overline{z}}+\sum_{j,k}\Gamma_{jk}^i\overline{\phi_j}\phi_k\Big\}\frac{\partial}{\partial
x_i}=0.
\end{aligned}\end{equation}
\end{rem}

In general it is quite difficult to produce functions $\phi_i$ with
the above properties since the holomorphicity condition is given by
partial differential equations with {\em nonconstant coefficients}.
If $M$ is a Lie group equipped with a left-invariant metric $g$ and
$\{E_i\}$ are orthonormal left-invariant vector fields, we can write
$$\phi = \sum_i \phi_i\frac{\partial}{\partial x_i} = \sum_i \psi_i E_i,\qquad \psi_i:
\Omega \rightarrow \C,$$ with $\phi_i = \sum_{i,j} A_{ij}\psi_j$ and
$A = (A_{ij})$ being an invertible matrix, with function entries
$A_{ij}$. In this case the Weierstrass formula becomes:
\begin{thm}{\em(see \cite{MMP} for a proof)}\label{wlie}
Given functions $\psi_i: \Omega \rightarrow  \C$ such that:
\begin{enumerate}
\item $\sum_i |\psi_i|^2 \neq 0$,
\item $\sum_{i}\psi_i^2 = 0$,
\item $\displaystyle\frac{\partial {\psi}_i}{\partial \overline z} +
\sum_{j,k}L_{jk}^i \overline{\psi_j}\psi_k=0,$
\end{enumerate}
where $L^i_{jk} :=g(\nabla_{E_j} E_k,E_i)$, then the map:
$$f: \Omega \rightarrow M,\qquad f_i(z) = 2\,\Re\bigg(\int_{z_0}^z \sum_j A_{ij}\psi_j\,
dz\bigg),$$
defines a conformal minimal immersion.
\end{thm}

 The advantage of having  partial differential equations with constant coefficients is not really a great gain,
 in principle,
since we still have to compute the integrand $A_{ij}\psi_j$ along
the solutions. However, in certain cases, as for example the
hyperbolic space, the Heisenberg group and  $\h^2\times\r$, this
problem may be overcome by ad hoc arguments, as shown (for example)
in \cite{MMP}.

\section{The Bj\"{o}rling problem for three-dimensional Lie groups}

In this section we will suppose that $M$ is a three-dimensional Lie
group endowed with a left-invariant Riemannian metric $g$. Let
$\beta:I\subseteq \r\rightarrow M$ be a regular analytic curve in
$M$ and $V:I\rightarrow TM$ a unitary real analytic vector field
along $\beta$, such that $g(\dot\beta,V)\equiv 0$. The Bj\"{o}rling
problem is the following:

\mbox{}

\noindent {\em Determine a minimal surface $f: I\times
(-\epsilon,\epsilon) = \Omega \subseteq \C\rightarrow M$, such that:
\begin{itemize}
  \item $f(u,0)=\beta (u)$,
  \item $N(u,0)=V(u)$,
\end{itemize}
for all $u\in I$, where $N: \Omega \rightarrow TM$ is the Gauss
map of the surface.}

\mbox{}

We observe that if $\beta$ is parameterized by arc-length and
$\ddot{\beta}:=\nabla_{\dot{\beta}} {\dot{\beta}}$, we have that
$V=\|\ddot{\beta}\|^{-1}\ddot{\beta}$ is a unit vector field along
the curve such that $g(\dot\beta,V)\equiv 0$. Then the Bj\"{o}rling
problem is a generalization of the problem of finding a minimal
surface which contains a given curve as a geodesic.

\begin{thm}\label{bjor} The Bj\"{o}rling
problem has a unique solution\footnote{Unique up to fixing the
domain.}.
\end{thm}
\begin{proof} In order to prove the theorem, we must analyse Theorem~\ref{wlie}
carefully. In this theorem we have essentially four conditions on
the three functions $\psi_i$ (the first condition is ``generically
satisfied''). We will start showing that these conditions are
dependent.
\begin{lem} \label{dipendenti}
Let $\psi_i: \Omega \subseteq\C\to\C$, $i = 1, 2,$ be two
differentiable functions and $\psi_3^2 = -\psi_1^2 - \psi_2^2$. We
suppose that $\psi_i$, $i = 1, 2$, satisfy the two first equations
of the third item of Theorem~\ref{wlie}. Then $\psi_3$ satisfies the
third equation.
\end{lem}
\begin{proof}
Deriving with respect to $\bar{z}$ the equation
$$-\psi_3^2=(\psi_1^2+\psi_2^2),$$ and using the fact that  the two first
equations of  the third item of Theorem~\ref{wlie} are satisfied, we
have:
\begin{align*} -\psi_3\,\frac{\partial\psi_3}{\partial \bar{z}}&=
\psi_1\,\frac{\partial\psi_1}{\partial
\bar{z}}+\psi_2\,\frac{\partial\psi_2}{\partial
\bar{z}}\\&=-\sum_{j,k=1}^3(L_{jk}^1\,\psi_1+L_{jk}^2\,\psi_2)\,\bar{\psi}_j\,\psi_k.
\end{align*}
Therefore, to prove the lemma it suffices to show that
$$\sum_{j,k=1}^3(L_{jk}^1\,\psi_1+L_{jk}^2\,\psi_2+L_{jk}^3\,\psi_3)\bar{\psi}_j\,
\psi_k=0.$$
Writing the above sum as:
$$\sum_{j,k=1}^3L_{jk}^k\,\bar\psi_j\,\psi_k^2+\sum_{\tiny\begin{array}{c}
  j,k,l=1 \\
  k<l
\end{array}}^3(L_{jk}^l+L_{jl}^k)\,\bar{\psi}_j\,\psi_k\,\psi_l,$$
and using the relation $L_{jk}^l+L_{jl}^k=0$, where
$j,k,l\in\{1,2,3\},$ we conclude the proof.
\end{proof}

We go back now to the proof of Theorem~\ref{bjor}. Consider the
system:
\begin{equation}\label{ya}
\left\{ \begin{aligned}
\frac{\partial \psi_1}{\partial \overline{z}} +
\sum_{j,k =1}^3L^1_{jk}\overline{\psi_j}\psi_k&= 0,\\
\frac{\partial \psi_2}{\partial \overline{z}} +
\sum_{j,k =1}^3L^2_{jk}\overline{\psi_j}\psi_k&= 0,\\
\end{aligned}
\right.
\end{equation}
where $\psi_i: \Omega \rightarrow \C$ and $\psi_3^2 = -\psi_1^2 -
\psi_2^2$.

Since this system is of Cauchy-Kovalevskaya type (see \cite{P} for a
proof of the Cauchy-Kovalevskaya Theorem), fixing the initial datas
$\psi_i(u,0)$, $ i = 1, 2$, it has, locally, a unique solution. This
solution gives, via Theorem~\ref{wlie} and Lemma~\ref{dipendenti}, a
minimal surface. Thus we must find initial conditions so that this
surface has the required properties. Observe that, if $f$ is a
solution of the Bj\"{o}rling problem, we have:
\begin{equation}\label{ci}\phi(u,0) := \frac{1}{2}\Big(\frac{\partial
f}{\partial u} - i\frac{\partial f}{\partial v}\Big)(u,0) =
\frac{1}{2}\big(\dot{\beta}(u) +i\,\dot{\beta}(u)\wedge
V(u)\big).\end{equation} Therefore, the initial data for the system
is:
\begin{equation}\label{ci1}\psi(u,0) = A^{-1}(\beta(u))\,\phi(u,0).
\end{equation}

Note that the initial condition implies
$\displaystyle{\frac{\partial f}{\partial u}(u,0) =
\dot{\beta}(u)}$. Hence, up to a constant determined by the constant
of integration in Theorem~\ref{wlie}, we have $f(u,0) = \beta(u)$.
Also, the initial condition forces the choice of one of the
determinations of $\psi_3^2 = -\psi_1^2 - \psi_2^2$.

Up to now we have proved the existence of a local solution to the
problem. Using compactness of $I$ and local uniqueness, we have
existence and uniqueness of the solution when $\beta(I)$ is
contained in a coordinate neighborhood, for $\epsilon$ sufficiently
small. Covering $I$ with a finite number of inverse images, via
$\beta$, of coordinate neighborhoods and using (again) the
uniqueness of the local problem, the result is proved for the
general case.
\end{proof}

\section{Examples in the space $\h^2 \times \r$}

 Let $\h^2$ be the hyperbolic plane
$\{(x,y)\in\r^2\;:\; y>0\}$ endowed with the metric, of constant
Gauss curvature $-1$, given by $g_{\tiny \h}=(dx^2+dy^2)/y^2$. The
hyperbolic plane $\h^2$, with the group structure derived by the
composition of proper affine maps, is a Lie group and the metric
$g_{\tiny \h}$ is left-invariant. Then the product space $\h^2\times
\r$ is a Lie group with the product structure given by
$$(x,y,z)*(x',y',z')=(x' y+x,y y',z+z')$$
 and the
product metric $ g=g_{\tiny \h}+dz^2 $ is left-invariant. The Lie
algebra of the infinitesimal isometries of $(\h^2\times\r,g)$ admits
the following bases of Killing vector fields
\begin{align*}
X_1&=\frac{(x^2-y^2)}{2}\frac{\partial}{\partial
x}+xy\frac{\partial}{\partial
y},\\
X_2&=\frac{\partial}{\partial x},\\
X_3&=x\frac{\partial}{\partial x}+y\frac{\partial}{\partial y},\\
X_4&=\frac{\partial}{\partial z}.
\end{align*}
 With respect to the metric $g$ an
orthonormal basis of left-invariant vector fields is:
$$ E_{1} = y\,\frac{\partial}{\partial
x},\quad  E_{2} = y\,\frac{\partial}{\partial y},\quad  E_{3} =
\frac{\partial}{\partial z}.$$ Also, the matrix $A$ is given by:
$$ A=\left(
\begin{array}{ccc} y&0&0\\ 0&y&0\\ 0&0&1
\end{array} \right),
$$ and the non zero $L_{ij}^k$ are $L_{12}^1=-1$ and $L_{11}^2=1$.
Consequently system~\eqref{ya} becomes:
\begin{equation}\label{ya1}
\left\{ \begin{aligned} \frac{\partial \psi_1}{\partial
\overline{z}} -\overline{\psi_1} \psi_2&= 0,\\
\frac{\partial \psi_2}{\partial \overline{z}} +
|\psi_1|^2&= 0.\\
\end{aligned}
\right.
\end{equation}

\begin{ex}[The horizontal plane $z=c$]
First of all, we consider the curve $$\beta (u)=(\cos u,\sin u,
c),\qquad u\in (0,\pi),\quad c\in\r,$$ and the unit vector field
$V(u)=-E_3$. As $\dot{\beta}(u)=-\sin u\,E_1+\cos u\, E_2$, it
results that $g(\dot{\beta},V)\equiv 0$ and, also,  using the
equations {\eqref{ci}} and {\eqref{ci1}}, we have that the initial
data for the system~{\eqref{ya1}} is: $$\psi (u,0)=\Big(-\frac{(\sin
u+ i\,\cos u) }{2\,\sin u},\frac{(\cos u-i\,\sin u)}{2\,\sin u},0
\Big).$$ Thus, it follows that $\psi (u,v)=\psi(u,0)$ and,
integrating, we obtain the conformal immersion of the totally
geodesic plane $z=c$ given by
$$f(u,v)=(e^v\,\cos u, e^v\,\sin u,c).$$
\end{ex}

\begin{ex}[The helicoid]
Consider the curve $\beta (u)=(0,1,2 u)$ and the unit vector field
$V(u)=\cos (2 u)\,E_1+\sin (2 u)\, E_2$. As $\dot{\beta}(u)=2\,E_3$,
it results that $g(\dot{\beta},V)\equiv 0$. Also, using the
equations {\eqref{ci}} and {\eqref{ci1}}, we have that the initial
data for the system~{\eqref{ya1}} is:
$$\psi (u,0)=(-i\,\sin (2u),i\,\cos (2 u),1)$$ and moreover $\psi_3 = (-\psi_1^2 -
\psi_2^2)^{1/2}$. Consequently, the solution is given by
\begin{equation}\begin{aligned}\psi_1(u,v)&=\frac{2\,i\,\sin(2 u) - 2\,(\cos(2 u) + \sin(2
v))\,\tan (2 v)}{ -2+\sin (2u -2v) - \sin (2u + 2v)},\\
\psi_2(u,v)&=\frac{\sec(2 v)\,[2\, i\,(\cos (2 u) + \sin (2 v)) +
\sin (2 u)\sin (4 v)]}{2 - \sin(2u - 2v) + \sin(2u + 2v)},\\
\psi_3(u,v)&=1.
\end{aligned}
\end{equation}
After integration, we have the immersion of the minimal helicoid in
$\h^2\times\r$ described in \cite{MMP}, given by:
\begin{equation}\nonumber\begin{aligned}f_1(u,v)&=\frac{2\,\sin(2 u)\,\sin(2 v)}{ 2-\sin (2u -2v) + \sin (2u + 2v)},\\
f_2(u,v)&=\frac{2\,\cos (2 v)}{2 - \sin(2u - 2v) + \sin(2u + 2v)},\\
f_3(u,v)&=2\,u.
\end{aligned}
\end{equation}
\end{ex}

\section{Examples in the Heisenberg group $\h_3$}
We now consider the Heisenberg group $$\h_3=\Bigg\{\left(
\begin{array}{ccc}
1&x&z+\frac{1}{2}x y\\ 0&1&y\\ 0&0&1
\end{array}
\right),\;x,y,z\in\r\Bigg\},
$$
equipped with the left-invariant metric given by
$$g=dx^2+dy^2+\Big(\frac{1}{2}ydx-\frac{1}{2}xdy+dz\Big)^2.$$ An orthonormal basis of left-invariant vector fields
is given by:
$$
E_{1} = \frac{\partial}{\partial x} - \frac{y}{2}
\frac{\partial}{\partial z},\qquad E_{2} = \frac{\partial}{\partial
y}+ \frac{x}{2} \frac{\partial}{\partial z},\qquad E_{3} =
\frac{\partial}{\partial z}.
$$
The matrix $A$ takes the form $$A= \left(
\begin{array}{ccc}
1&0&0\\ 0&1&0\\ -\frac{y}{2}&\frac{x}{2}&1
\end{array}
\right)
$$
and the non zero $L_{ij}^k$ are:
$$\begin{aligned}&L_{12}^3=\frac{1}{2},\qquad
L_{21}^3=-\frac{1}{2},\\&L_{13}^2=-\frac{1}{2},\quad
L_{31}^2=-\frac{1}{2},\\&L_{23}^1=\frac{1}{2},\qquad
L_{32}^1=\frac{1}{2}.
\end{aligned}$$
Thus, the system~\eqref{ya} becomes:
\begin{equation}\label{ya2}
\left\{ \begin{aligned} \frac{\partial \psi_1}{\partial
\overline{z}} +\Re(\psi_2\,\overline{\psi_3})&= 0,\\
\frac{\partial \psi_2}{\partial \overline{z}}-\Re(\psi_1\,\overline{\psi_3})&= 0.\\
\end{aligned}
\right.
\end{equation}

\begin{ex}[Helicoids]
We consider $$\beta (u)=(\rho (u),0,b)\qquad\text{and}\qquad
V(u)=\frac{(\rho^2-2c)}{2 \rho'}\, E_2+\frac{\rho }{\rho'}\,E_3,$$
where $b,c\in\r$ and the real-valued function $\rho=\rho (u)$
satisfies:
$$\sqrt{(\rho')^2-\rho^2}=\rho^2/2-c.$$ Since
$\dot{\beta}(u)=\rho'(u)E_1$, we obtain that $g(\dot{\beta},V)\equiv
0$. Using {\eqref{ci}} and {\eqref{ci1}}, we have that the initial
data for the system~{\eqref{ya2}} is
$$\psi(u,0)=\frac{1}{2}\big(\rho'(u),-i\, \rho(u),i\, (\rho(u)^2-2c)/2\big).$$ Consequently $\psi_3 = (-\psi_1^2 -
\psi_2^2)^{1/2}$ and $$\psi (u,v)=\frac{1}{2}\big(\rho'(u)\,\cos v+i
\,\rho(u)\, \sin v, \rho'(u)\,\sin v-i\, \rho(u)\, \cos
v,i\,(\rho(u)^2-2 c)/2\big).$$ Using the fact that:
$$\phi_1=\psi_1,\qquad
\phi_2=\psi_2,\qquad\phi_3=-\frac{y}{2}\,\psi_1+\frac{x}{2}\,\psi_2+\psi_3,$$
and:
\begin{equation}\label{lai}
\left\{ \begin{aligned}f_1(z) &=2\,\Re\int_{z_0}^z \psi_1\,dz,\\
f_2(z)&=2\,\Re\int_{z_0}^z \psi_2\,dz,\\
f_3(z)&=2\,\Re\int_{z_0}^z
\Big(\psi_3-\frac{f_2}{2}\,\psi_1+\frac{f_1}{2}\,\psi_2\Big)\,dz,
\end{aligned} \right.\end{equation}
we obtain the minimal immersion:
$$f(u,v)=(\rho (u)\,\cos v,\rho (u)\,\sin v,c\,v+b).$$ Therefore, if $c\neq 0$
we have the parametrization of a helicoid, while, if $c=0$ we
obtain the horizontal plane $z=b$.
 \end{ex}

\begin{ex}[Catenoid-type surface]
We consider the curve in $\h_3$ given by: $$\beta (u)=(g\,\cos l,
\,g\,\sin l,\tilde h),$$ where $g=g(u)$, $l=l(u)$ and $\tilde
h=\tilde h (u)$ are real-valued functions such that
$$g'^2=\frac{g^2(g^4-16)-4}{g^2-4}$$ and
$${\tilde h}'=h=\sqrt{\frac{g^2+4}{g^2-4}},\qquad l'=\frac{2h}{g^2+4},$$ with
$g^2>4$. Let
$$V(u)=-\frac{(gg'\sin l+2h\cos l)}{g\,(g^2+4)}\,E_1+\frac{(gg'\cos l-2h\sin l)}{g\,(g^2+4)}\,E_2
+\frac{2g'}{g\,(g^2+4)}\,E_3$$ be a unitary vector field. As
$$\dot{\beta}(u)=(g'\cos l-g\,l'\sin l)\,E_1+(g'\sin l+g\, l'\cos l)\,E_2+2l'\, E_3,$$ it is easy to check that $g(\dot{\beta},V)\equiv
0$ and $$\dot{\beta}\wedge V=2g\sin l\,E_1-2 g\cos
l\,E_2+g^2\,E_3.$$ Therefore, using {\eqref{ci}} and {\eqref{ci1}},
it follows that the initial data for system~{\eqref{ya2}} is
$$\psi(u,0)=\frac{1}{2}(g'\cos l-gl'\sin l+2 i g\sin l,g'\sin l+g l'\cos l-2 i g \cos l, 2l'+i g^2).$$
Hence,
\begin{equation}\nonumber
\begin{aligned}\psi_1 (u,v)&=\frac{1}{2}(g'\cos (l+2v)-g l' \sin (l+2v)+2ig\sin
(l+2v)),\\\psi_2 (u,v)&=\frac{1}{2}(g'\sin (l+2v)+g l' \cos
(l+2v)-2ig\cos (l+2v)),\\ \psi_3(u,v)&=\frac{1}{2}(2l'+ig^2).
\end{aligned}\end{equation}
After integration, we have the catenoid-type minimal surface given
by $$f(u,v)=(g\, \cos (l+2v), g\,\sin (l+2v),\widetilde{h}).$$
\end{ex}

\section{Final comments}

If the ambient space is $\r^3$ with the flat metric, the solution of
the Bj\"{o}rling problem can be given by an explicit formula. This
fact has been used to prove a {\em reflection principle} and a nice
application of this is the characterization of the helicoid as the
unique ruled minimal surface, besides the plane (see \cite{H},
\cite{S}).

In our case, the partial differential equations involved are more
complicated than the Cauchy-Riemann ones, and we were not able to
find an ``explicit'' formula for the solution of the Bj\"{o}rling
problem. It is not clear that a reflection principle holds for a generic three-dimensional Lie group.
For example, in the case of the Heisenberg group
$\h_3$, there are many minimal surfaces ruled by translated of
$1$-parameter subgroups. Such surfaces were classified in \cite{B}
and \cite{F}. If we consider the developable ones, i.e. the ones
with Gauss map of rank $1$, then they are (up to isometries of the
ambient space) the graphs of the following functions:
\begin{equation}\nonumber
f(x,y)=\left\{
\begin{array}{l}
{\displaystyle\frac{xy}{2}+k\,[\ln(y+\sqrt{1+y^2})+y\,\sqrt{1+y^2}]}\\
\text{or}\\
{\displaystyle 2 k y-\frac{xy}{2}},\qquad k\in\r.\\
\end{array}
\right.
\end{equation}

Since these are complete graphs, the Bernstein Theorem does not
hold. A classification of complete minimal graphs has been recently
given in \cite{FM}, in terms of the generalized Hopf differentials
introduced in \cite{AR}.

We could also consider minimal surfaces of a Lie group, in
particular of $\h_3$, ruled by geodesics. Very little is known for
such surfaces and the classification problem seems to be more
difficult. However, our feeling is that this is the right context
for a reflection principle.

\end{document}